\documentclass[12pt]{amsart}

\usepackage{graphicx}
\usepackage{amsthm}
\usepackage{amsmath, amssymb}
\swapnumbers

\theoremstyle{plain}
\newtheorem{Thm}{Theorem}

\newtheorem{Coro}[Thm]{Corollary}
\newtheorem{Lem}[Thm]{Lemma}
\newtheorem{Claim}{Claim}

\theoremstyle{definition}
\newtheorem{Def}[Thm]{Definition}

\begin{document}

\title{Surface bundles with genus two Heegaard splittings}
\author{Jesse Johnson}
\address{\hskip-\parindent
        Department of Mathematics\\
        University of California\\
        Davis, CA 95616\\
        USA}
\email{jjohnson@math.ucdavis.edu}
\subjclass{Primary 57M}
\keywords{Heegaard splitting, surface bundle, Rubinstein-Scharlemann graphic}

\thanks{Research supported by NSF VIGRE grant 0135345}

\begin{abstract}
It is known that there are surface bundles of arbitrarily high genus which have genus two Heegaard splittings.  The simplest examples are Seifert fibered spaces with the sphere as a base space, three exceptional fibers and which allow horizontal surfaces.  We characterize the monodromy maps of all surface bundles with genus two Heegaard splittings and show that each is the result of integral Dehn surgery in one of these Seifert fibered spaces along loops where the Heegaard surface intersects a horizontal surface.  (This type of surgery preserves both the bundle structure and the Heegaard splitting.)
\end{abstract}

\maketitle

\section{Introduction}

A \textit{handlebody} is a 3-manifold (with boundary) that is homeomorphic to a closed regular neighborhood of a connected graph embedded in a compact, orientable 3-manifold.  Given a compact, closed, connected and orientable 3-manifold, $M$, let  $\Sigma$ be a surface embedded in $M$ and let $H_1$ and $H_2$ be handlebodies embedded in $M$.  The triple $(\Sigma, H_1, H_2)$ is a \textit{Heegaard splitting} if $H_1 \cup H_2 = M$ and $H_1 \cap H_2 = \Sigma = \partial H_1 = \partial H_2$.  The \textit{genus} of a Heegaard splitting is the genus of $\Sigma$.

An orientable Seifert fibered space with base the 2-sphere and three or fewer exceptional fibers is often called a \textit{small} Seifert Fibered space because it contains no closed incompressible surfaces.  Given a small Seifert fibered space $M$, one can construct a genus two Heegaard splitting for M as follows:  Let $K$ be the union of one exceptional fiber and a horizontal arc that passes between the other two exceptional fibers.  Let $H_1$ be the closure of a regular neighborhood of $K$.  The closure of the complement of $H_1$ is a second handlebody, $H_2$.  Letting $\Sigma = \partial H_1$ then $(\Sigma, H_1, H_2)$ is a Heegaard splitting for $M$.  See~\cite{lum:niels} for a more detailed description of this construction.

The construction works for any choice of coefficients for the three exceptional fibers.  By choosing the coefficients carefully, one can ensure that there is a horizontal incompressible surface $S \subset M$ of arbitrary genus.  In this case, $M$ also has the structure of a surface bundle over $S$.  That is, there's a smooth map $\pi : M \rightarrow S^1$ such that the preimage of any point in $S^1$ is a closed surface in $M$ isotopic to $S$.  We can reconstruct $M$ from $S \times [0,1]$ by gluing  $S \times \{0\}$ to $S \times \{1\}$ as defined by some map $\phi : S \rightarrow S$ called the \textit{monodromy} of $M$.

Assume $S$ and $\Sigma$ are transverse and let $\ell$ be a simple closed curve in $\Sigma \cap S$.  For a neighborhood $N$ of $\ell$, let $\lambda$ be a longitude defined by $\partial N \cap S$.  This $\lambda$ is also isotopic to a component of $\partial N \cap \Sigma$.  Let $M'$ be the result of $n$ surgery on $\ell$.  In other words, we remove $N$ then replace it, adding $n$ Dehn twists along $\lambda$ to the gluing map.  

This construction is equivalent to adding $n$ Dehn twists to the monodromy map of the surface bundle along $\ell \subset S$.  It is also equivalent to adding $n$ Dehn twists to the gluing map from $\partial H_1$ to $\partial H_2$ along $\ell \subset \Sigma$.  This is because $\lambda$ can be made to lie in $S$ or $\Sigma$.  As a result, $M'$ is a new surface bundle over $S$ and $M'$ allows a genus two Heegaard splitting.  (However, $M'$ may no longer be a Seifert fibered space.)

A second way to construct surface bundles with genus two Heegaard splittings is as follows:  Given a compact, closed, orientable surface $S$ of genus $g$, let $G \subset S$ be a graph consisting of a single vertex and $2g$ edges, such that the complement of $G$ is a single open disk.  We will call $G$ a \textit{one-vertex spine} of $S$.  An automorphism  $\phi : S \rightarrow S$ is a \textit{cyclic permutation of $G$} if for some labeling $e_1,\dots,e_{2g}$ of the edges of $G$, $\phi$ sends each edge $e_i$ onto $e_{i+1}$ and sends $e_{2g}$ onto $e_1$.  

Let $N_1$ be a regular neighborhood of $e_1$ in $S$.  Let $N_2$ be a regular neighborhood of $e_1 \cup e_2$ containing $N_1$ in its interior and similarly, for each $i < 2g$, let $N_i$ be a regular neighborhood of $e_1,\dots, e_i$ containing $N_{i-1}$ in its interior.

We will prove the following theorem:

\begin{Thm}
\label{mainthm}
Let $M$ be an orientable surface bundle over a compact, connected, closed, orientable surface $S$ with monodromy $\phi : S \rightarrow S$.  Then $M$ allows a genus two Heegaard splitting if and only if for some one-vertex spine $G$ of $S$, $\phi$ is isotopic to the composition of a cyclic permutation of $G$ and a number of Dehn twists along the boundaries of $N_1,\dots,N_{2g-1}$.
\end{Thm}

\begin{Coro}
\label{maincoro}
The examples constructed above from small Seifert fibered spaces by performing Dehn surgery along the intersection of the Heegaard surface and a horizontal surface are the only surface bundles with genus two Heegaard splittings.
\end{Coro}

Scharlemann and Cooper~\cite{sch:solv} have classified Heegaard splittings for torus bundles. The case of Theorem~\ref{mainthm} when $S$ is a torus follows from their results.  Our proof of Theorem~\ref{mainthm} will thus be restricted to the case when $S$ has genus two or greater.  Schleimer and Bachman~\cite{bsc:bndls} have shown that the genus of a strongly irreducible Heegaard splitting of a surface bundle is bounded below based on how the monodromy map acts on the curve complex.  The proof here uses a similar technique, with a deeper analysis of the genus two case.

The proof of Theorem~\ref{mainthm} is broken into two distinct parts.  In the first half of the paper, Sections~\ref{constrsect} through~\ref{mainthmsect}, we assume the Heegaard surface has been isotoped so as to intersect the level surfaces of the bundle structure in a ``nice'' way, then use cut-and-paste style arguments to characterize the monodromy of the surface bundle.  In the second half of the paper, Sections~\ref{positionsect} through~\ref{fairsect}, we use the double sweep-out method of Rubinstein and Scharlemann~\cite{rub:compar} to show that a genus two Heegaard surface can always be isotoped to meet the criteria needed for the first half.

\section{Constructing a genus two Splitting}
\label{constrsect}

In this section we will show that a surface bundle $M$ with leaf $S$ and monodromy, $\phi$, as described in Theorem~\ref{mainthm}, allows a genus two Heegaard splitting.  We will then prove that Corollary~\ref{maincoro} follows from Theorem~\ref{mainthm}.  Let $S$ be a compact, closed, connected, orientable surface, $G$ a one-vertex spine for $S$ and $\phi$ an automorphism of $S$ as described in Theorem~\ref{mainthm}.  

Let $v$ be the vertex and $e_1,\dots,e_{2g}$ the edges of $G$.  Write $\phi = \phi_{2g-1} \circ \phi_{2g-2} \circ \dots \circ \phi_{0}$ where $\phi_0$ is a cyclic permutation of $G$ and for each $i < 2g$, $\phi_i$ is some composition of Dehn twists along the boundary loops of a regular neighborhood of  the edges $e_1 \cup \dots \cup e_i$ of $G$.  These loops are pairwise disjoint so the order of the dehn twists does not matter.

We will construct a sequence of surface bundles $M_0,\dots,M_{2g-1}$  and show that each has a genus two Heegaard splitting.  In particular, each $M_i$ will be the result of a Dehn  surgery on $M_{i-1}$ along simple closed curves of intersection between the Heegaard surface and a leaf.  The final manifold, $M_{2g-1}$, will be homeomorphic to $M$.

Let $U$ be the universal cover of $S$.  The spine of $S$ lifts to a graph which defines a tiling of the plane $U$ by $2g$-gons.  Let $P$ be one of these polygons.  This $P$ is a fundamental domain for the covering of $S$ by $U$ and we can reconstruct $S$ by identifying pairs of edges of $P$.  The automorphism $\phi_0$ sends the spine $G$ of $S$ onto itself, so it defines an automorphism of $P$.  

If we identify $P$ with a regular polygon in the plane, centered at the origin, then any automorphism of $P$ is isotopic to a rotation of angle $k\pi/g$ for some $k$.  For any $i$, $j$, some power of $\phi_0$ will send $e_i$ to $e_j$, so for any edges $a$, $b$ of $P$, the action of some power of $\phi_0$ will send $a$ to either $b$ or the edge to which $b$ is identified.  

The action of $\phi_0$ on $P$ sends pairs of edges identified to pairs of identified edges and some power of $\phi_0$ sends each pair to each pair.  This is only possible if each edge is identified with the edge opposite it in $P$.  Consider an edge $a$ of $P$ and let $b$ be an adjacent edge.  Some power $\phi_0^l$ sends $a$ to either $b$ or the edge opposite $b$.  In either case, the order of $\phi_0^l$ is $2g$, so the order of $\phi_0$ is at least $2g$.  In fact, this is the most the order can be, so $\phi_0$ has order exactly $2g$.  The map $\phi_0^g$ sends $e_1$ onto itself, but with opposite orientation.

The action of $\phi_0$ fixes the center of $P$ and sends any other point in the interior back to itself after $2g$ iterations.  This implies that in $S$ there is a single point in the complement of the spine $G$ that is fixed by $\phi_0$ and every other point has order $2g$.  The vertex of $G$ is fixed by $\phi_0$.  The power of $\phi_0$ that sends an edge of $P$ onto its opposite sends the corresponding edge of $G$ onto itself, but with opposite orientation.  Thus the center of each edge has order $g$ and each non-center point in an edge has order $2g$.

Let $M_0$ be the result of gluing the boundary surface $S \times \{0\}$ of $S \times [0,1]$ to $S \times \{1\}$ according to the map $\phi_0$.  In other words, identify $(x, 0)$ with $(\phi_0(x),1)$ for each $x \in S$.

For each point $p \in S$, with orbit $A = \{\phi_0^i(p)\} \subset S$, the image in $M_0$ of $A \times [0,1]$ is a simple curve.  Because each point has finite orbit, each such curve is a closed loop and this family of loops defines a Seifert fibered structure on $M_0$.  The Seifert structure has three singular fibers:  The center of $P$ and the vertex of $G$ are each fixed by $\phi_0$ and each forms one singular fiber.  The set of centers of the edges of $G$ is a single orbit of $\phi_0$ and forms the third singular fiber.

Let $l_1$ be the image in $M_0$ of the loop $e_1 \times \{0\}$ and let $l_2$ be the inclusion in $M_0$ of the arc $\{v\} \times [0,1]$.  Then $l_2$ is a loop in $M_0$ and $l_1 \cup l_2$ is a graph with two edges and a single valence four vertex.  Let $N$ be an open regular neighborhood in $M_0$ of $l_1 \cup l_2$ and let $H^0_2$ be its complement.  Let $H^0_1$ be the closure of $N$ and define $\Sigma^0 = \partial H^0_1$.   We will show that $H^0_2$ is a handlebody, implying that the triple $(\Sigma^0, H^0_1, H^0_2)$ is a Heegaard splitting of $M_0$.

For each edge $e_i$ of $G$, the inclusion of the set $(e_i \cup v) \times [0,1]$ into $M_0$ is an annulus $T_i$.  The complement $D_i = T_i \setminus N$ is a disk and the intersection $D_i \cap D_{j}$ is an arc when $j = i \pm 1$ and is empty otherwise.  Thus $D = \bigcup D_i$ is a disk which is properly embedded in $H_2$.   The complement in $M_0$ of $H_1^0 \cup D$ is a subset of the open solid torus $(S \setminus G) \times S^1$ and is in fact homeomorphic to an open solid torus.  The set $H^0_2$ is the result of attaching a one-handle to this solid torus, so $H^0_2$ is a handlebody and $M_0$ allows a genus two Heegaard splitting $(\Sigma^0, H^0_1, H^0_2)$.

Let $N_1$ be a closed regular neighborhood of $D_1$ in $H^0_2$.  The disk $D_1$ intersects $\Sigma^0$ in a single arc so $N_1$ intersects $\Sigma^0$ in a disk.  Because $N_1$ is a ball and $N_1 \cap H^0_1$ is a disk, the union $H_1 \cup N_1$ is a handlebody and the closure of its complement is a second handlebody.  The intersection of $H^0_1 \cup N_1$ and $S \times \{\frac{1}{2}\}$ is a regular neighborhood of $e_1$ so Dehn surgery along the boundary loops will produce a manifold $M_1$ with monodromy $\phi_1 \circ \phi_0$.  

The inclusion of $H^0_1 \cup N_1$ into $M_1$ is a handlebody $H^1_1$ and the inclusion of the closure of $H^0_2 \setminus N_1$ into $M_1$ is a second handlebody $H^1_2$.  Thus $M_1$ allows a genus two Heegaard splitting $(\Sigma^1, H^1_1, H^1_2)$.  

The Dehn surgery does not affect the topology of $H^0_2$ so the inclusion of the disks $D_2 \setminus N_1$ and $D_3,\dots, D_{2g-1}$ into $M_1$ is a new collection of disks whose union is a meridian disk for $H^1_2$.  A regular neighborhood $N_2$ of $D_2$ will intersect $\Sigma^1$ in a disk and the intersection $(H^1_1 \cup N_2) \cap S \times \{\frac{1}{2}\}$ is a regular neighborhood of $e_1 \cup e_2$.  Thus we can repeat the construction to produce a manifold $M_2$ with monodromy $\phi_2 \circ \phi_1 \circ \phi_0$.  

By repeating this construction further, we eventually produce a manifold $M_{2g-1}$ with a genus-2 Heegaard splitting and monodromy $\phi_{2g-1} \circ \phi_{2g-2} \circ \dots \circ \phi_{0}$.  This manifold has the same monodromy as $M$ so $M_{2g-1}$ is homeomorphic to $M$.

\begin{proof}[Proof of Corollary~\ref{maincoro}]
Let $M$ be a surface bundle with fiber $S$, which allows a genus two Heegaard splitting. By Theorem~\ref{mainthm}, the monodromy map is of the form described above.  In the construction, we began with a small Seifert Fibered space and constructed $M$ by a sequence of Dehn surgeries along loops of intersection between the Heegaard surface and the leaves of the bundle.  This proves the corollary.
\end{proof}

\section{Outline of the Converse}
\label{outlinesect}

Let $M$ be a surface bundle and $K \subset M$ a graph.  An edge of $K$ is \textit{vertical} if it is transverse to every fiber and \textit{horizontal} if it is contained in a fiber.  In the construction in Section~\ref{constrsect}, the properties of the monodromy map are used to construct a spine for one of the handlebodies in $M$ consisting of a vertical edge and a horizontal edge.  The first step in the proof of the converse is to show that for any genus two Heegaard splitting, such a spine exists.

\begin{Lem}
\label{nicespine1lem}
Let $(\Sigma, H_1, H_2)$ be a genus two Heegaard splitting of $M$.  After an isotopy of $\Sigma$, there is a spine $K$ of $H_1$ such that $K$ consists of a single vertical edge and a single horizontal edge.
\end{Lem}

To prove this Lemma, we will consider the following construction:  Let $F$ be a surface with boundary.  Let $\alpha_0$ and $\alpha_1$ be properly embedded, essential arcs in $F$ (not necessarily disjoint). Let $N$ be an open regular neighborhood of $\alpha_0$ and the boundary components of $F$ containing the endpoints of $\alpha_0$.  Let $N'$ be a similar neighborhood for $\alpha_1$.  Let $\phi : (F \setminus N) \rightarrow (F \setminus N')$ be a homeomorphism.  

Define $N_0 = N \times [0,\frac{1}{4})$ and $N_1 = N' \times (\frac{3}{4},1]$ to be subsets of $F \times [0,1]$.  These are open neighborhoods of $\alpha_0 \times \{0\}$ and $\alpha_1 \times \{1\}$, respectively, such that $N_0 \cap (F \times \{0\}) = N \times \{0\}$ and $N_1 \cap (F \times \{1\}) = N' \times \{1\}$.  Define $X$ to be the quotient of $F \times [0,1] \setminus (N_0 \cup N_1)$ by the identification $(x,0) = (y,1)$ if $x \in F \setminus N$ and $\phi(x) = y$.

\begin{Def}
 Any manifold resulting from the above construction is an \textit{almost bundle}.
\end{Def}

Notice that the boundary of an almost bundle is a genus two surface.

\begin{Lem}
\label{spinelem}
If $X$ is an almost bundle which is homeomorphic to a genus two handlebody then there is a spine for $X$ consisting of a single vertical edge and a single horizontal edge, each of which forms a loop.
\end{Lem}

This Lemma will be proved in Sections~\ref{disksect} and \ref{structsect}.  To see its relevance to the main theorem, note the following lemma, which will be proved in Sections~\ref{positionsect} through \ref{fairsect}.

\begin{Lem}
\label{bundlehomeolem}
Let $(\Sigma, H_1, H_2)$ be a genus two Heegaard splitting for $M$.  After a suitable isotopy of $H_1$, there is an almost bundle $X$ and a homeomorphism $X \rightarrow H_1$ which takes leaves of the almost-bundle structure on $X$ into leaves of the bundle structure on $M$.
\end{Lem}

Because the proof of Lemma~\ref{bundlehomeolem} is self contained, it is presented at the end of the paper.

\section{Disks in Almost Bundles}
\label{disksect}

Let $F$ be a surface with boundary.  Let $\alpha_0$ and $\alpha_1$ be properly embedded, essential arcs in $F$ (not necessarily disjoint.)  As in the previous section, let $N$ be an open regular neighborhood in $F$ of $\alpha_0$ and the boundary components containing the endpoints of $\alpha_0$.  Let $N'$ be a regular neighborhood of $\alpha_1$ and the boundary components containing its endpoints.  

\begin{figure}[htb]
  \begin{center}
  \includegraphics[width=3.5in]{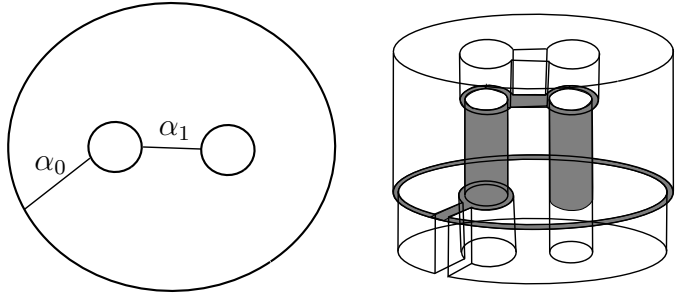}
  \put(-242,47){$\alpha_0$}
  \put(-195,60){$\alpha_1$}
  \caption{The surface on the left, with arcs $\alpha_0$ and $\alpha_1$, defines the manifold $Y$ on the right.  The shaded parts of $\partial Y$ plus an annulus around the outside form the set $Z$.}
  \label{xandyfig}
  \end{center}
\end{figure}

Let $N_0 = N \times [0,\frac{1}{4}) \subset F \times [0,1]$ and $N_1 = N' \times (\frac{3}{4},1] \subset F \times [0,1]$.  Define (as in Figure~\ref{xandyfig})
$$Y = F \times [0,1] \setminus (N_0 \cup N_1)$$ and
$$Z = (N \times \{1/4\}) \cup (N' \times \{3/4\}) \cup (\partial F \times [1/4,3/4]) \subset \partial Y.$$

\begin{Lem}
\label{techlem}
Assume there is a properly embedded disk $D \subset Y$ such that either $\partial D$ is contained in $Z$ or $\partial D$ is the union of an arc $a \subset Z$ and an arc $b$ disjoint from $Z$.  Then either

(1) The arc $\alpha_0$ can be isotoped disjoint from $\alpha_1$ or

(2) the boundary of $D$ bounds a disk in $\partial Y$.
\end{Lem}

\begin{proof}
If $\partial D$ is contained in $Z$ then let $a = \partial D$.  Otherwise, assume $\partial D$ is the union of an arc $a \subset Z$ and an arc $b$ with interior disjoint from $Z$.

Let $E$ be the properly embedded disk $(\alpha_0 \times [\frac{1}{4},1]) \cap Y$.   Assume $N'$ and $D$ have been isotoped so as to minimize $D \cap E$ while preserving the properties of $\partial D \cap Z$.  Assume we have isotoped $a$, so as to minimize the number of subarcs in each of the three pieces: $\partial F \times [\frac{1}{4},\frac{3}{4}]$, $N \times \{\frac{1}{4}\}$ and $N' \times \{\frac{3}{4}\}$.  Finally, assume that the intersecton $a \cap \partial E \cap (N \times \{\frac{1}{4}\})$ has been minimized and $a \cap E \cap \partial F \times [\frac{1}{4},\frac{3}{4}]$ has been minimized.

\begin{Claim}
\label{topnclaim}
Every arc of $a$ in $N' \times \{\frac{3}{4}\}$ that is not adjacent to $b$ is properly embedded and essential. If such an arc is disjoint from $E$ then $\alpha_0$ and $\alpha_1$ can be made disjoint.
\end{Claim}

\begin{proof}
If a subarc $\gamma$ of $a$ in $N' \times \{\frac{3}{4}\}$ is not adjacent to $b$ then the endpoints of $\gamma$ are in the interior of $a$, so $\gamma$ is properly embedded.  If $\gamma$ is properly embedded and essential in $F$ then $\gamma$ is parallel to $\alpha_1$ or cuts off an annulus containing $\alpha_1$ from $F$, as in Figure~\ref{essarcsfig}.  Thus if $\gamma$ is disjoint from $E$ then $\alpha_1$ is disjoint from $\alpha_0$.  A trivial arc in $N' \times \{\frac{3}{4}\}$ can be pushed into $\partial F \times [\frac{1}{4},\frac{3}{4}]$ so we can assume every subarc in $N' \times \{\frac{3}{4}\}$ is essential. 
\end{proof}
\begin{figure}[htb]
  \begin{center}
  \includegraphics[width=2.5in]{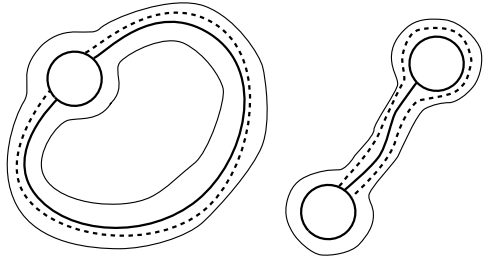}
  \caption{The dotted lines indicate the two types of essential subarcs of $a$ that can occur in $N' \times \{\frac{3}{4}\}$.  An arc parallel to $\alpha_1$ can occur whether the endpoints of $\alpha_1$ are in the same component of $\partial F$ or different component.  The second type of arc can only occur if the endpoints are in different components.}
  \label{essarcsfig}
  \end{center}
\end{figure}

\begin{Claim}
\label{bottomclaim}
Every arc of $a$ in $N \times \{\frac{3}{4}\}$ that is not adjacent to $b$ is properly embedded, essential and disjoint from $E$.
\end{Claim}

\begin{proof}
As with arcs in $N' \times \{\frac{3}{4}\}$, we can push trivial arcs out of $N \times \{\frac{1}{4}\}$, so that every arc in $N \times \{\frac{1}{4}\}$ is essential.  As in $N'$, each arc in $N$ is parallel to $\alpha_0$ or cuts off an annulus containing $\alpha_0$.  Thus the assumption that $a \cap \partial E \cap (N \times \{\frac{1}{4}\})$ is minimized implies that each arc is disjoint from $\alpha_0$ and therefore disjoint from $E$.
\end{proof}

\begin{Claim}
\label{annclaim}
Every arc of $a$ in $\partial F \times [\frac{1}{4},\frac{3}{4}]$ is essential and dsjoint from $E$.
\end{Claim}

\begin{proof}
The neighborhood $N$ contains a neighborhood of any boundary component that it touches, as does $N'$, so any trivial arc in $\partial F \times [\frac{1}{4},\frac{3}{4}]$ can be pushed into $N \times \{\frac{1}{4}\}$ or $N' \times \{\frac{3}{4}\}$.  Thus every arc in $\partial F \times [\frac{1}{4},\frac{3}{4}]$ must go from $N \times \{\frac{1}{4}\}$ to $N' \times \{\frac{3}{4}\}$.

Every arc of $(\partial F \times [\frac{1}{4},\frac{3}{4}]) \cap E$ is of the form $\{x\} \times [\frac{1}{4},\frac{3}{4}]$, i.e. ``vertical''.  Every arc of $a \cap (\partial F \times [\frac{1}{4},\frac{3}{4}]$ is essential and therefore isotopic to a vertical arc.  Moreover, we can choose an isotopy which fixes the endpoint of the arc in $\partial F \times \{\frac{1}{4}\}$ so that the induced isotopy of $a$ does not increase the number of points in $a \cap \partial E \cap (N \times \{\frac{1}{4}\}$.  Once all the arcs have been made vertical, they will be disjoint from $\partial E$.
\end{proof}

If some component of $D \cap E$ is a loop then this loop bounds a disk in $D$ and a disk in $E$.  We can compress $D$ along an innermost disk in $E$, forming a new disk with the same boundary as $D$ but with fewer intersections with $E$.  Because the number of intersections is minimized, $D \cap E$ must consist entirely of arcs.

\begin{Claim}
If $D$ is disjoint from $E$ then $\alpha_0$ and $\alpha_1$ can be made disjoint or $\partial D$ is trivial in $\partial Y$.
\end{Claim}

\begin{proof}
Assume $D$ is disjoint from $E$.  Then any arc of $a$ in $N' \times \{\frac{3}{4}\}$ is disjoint from $E$.  If such an arc exists then Claim~\ref{topnclaim} implies $\alpha_0$ and $\alpha_1$ can be made disjoint.   

Otherwise, assume $a$ is disjoint from $N' \times \{\frac{3}{4}\}$.   Any arc in $\partial F \times [0,1]$ with both endpoints in $N \times \{\frac{1}{4}\}$ is trivial so Claim~\ref{annclaim} implies that $a$ is contained in $N \times \{\frac{1}{4}\}$.  

If $a$ is all of $\partial D$ then the projection of $D$ into $F$ is an immersed disk $D' \subset F$ whose boundary coincides with $\partial D$.  A simple closed curve which bounds an immersed disk in a surface bounds an embedded disk so $\partial D$ bounds a disk in $\partial Y$.  If $b \neq \emptyset$ then the projections of $a$ and $b$ are disjoint since $a \subset N$ and $b \cap N = \emptyset$.  The projection of $a \cup b$ is simple in $F$ and bounds an immersed disk.  Once again, this loops also bounds an embedded disk in $\partial Y$.
\end{proof}

Assume that $D \cap E$ is not empty.  The arcs $D \cap E$ cut off at least two outermost disks from $D$.  The arc $b$ is disjoint from $\partial E$ so $b$ intersects at most one of the outermost disks.  Thus some arc $\zeta \subset (D \cap E)$ cuts off an outermost disk $D' \subset D$ such that $b$ is disjoint from $\partial D'$.  Let $\eta = \partial D' \cap \partial D$, so that $\partial D' = \zeta \cup \eta$ and $\eta$ is an arc in $Z$.

By Claims~\ref{bottomclaim} and~\ref{annclaim}, the arc $a$ can intersect $\partial E$ only in $N' \times \{\frac{3}{4}\}$.  The arc $\zeta$ is contained in $E$, so both endpoints of $\eta$ are in $N' \times \{\frac{3}{4}\}$ and $\eta$ is either contained in $N' \times \{\frac{3}{4}\}$ or goes from $N' \times \{\frac{3}{4}\}$ into $N \times \{\frac{1}{4}\}$ and then back into $N' \times \{3/4\}$.  If the arc leaves $N' \times \{\frac{3}{4}\}$ a second time then there is a subarc of $\partial D$ in $N' \times \{\frac{3}{4}\}$ which is properly embedded in $F$ and disjoint from $E$, so Claim~\ref{topnclaim} would imply that $\alpha_0$ and $\alpha_1$ are disjoint.

\begin{Claim}
\label{topclaim1}
If $\eta$ is contained in $N' \times \{\frac{3}{4}\}$ then $D \cap E$ can be reduced.
\end{Claim}

\begin{proof}
The projection of $\partial D'$ into $F$ is a simple closed curve which bounds an immersed disk (the projection of $D'$)  and therefore an embedded disk.  Sliding $\eta$ across this disk (and $N'$ with it) reduces the number of components of $D \cap E$.
\end{proof}

We assumed that $D \cap E$ is minimal, so we must have that $\eta$ begins in $N' \times \{\frac{3}{4}\}$, crosses $\partial F \times [\frac{1}{4},\frac{3}{4}]$ into $N \times \{\frac{1}{4}\}$,  then crosses $\partial F \times [\frac{1}{4},\frac{3}{4}]$ back into $N' \times \{\frac{3}{4}\}$.  If $\alpha_0$ and $\alpha_1$ are not disjoint then $\eta$ must end after it reenters $N' \times \{\frac{3}{4}\}$.  Let $l$ be the projection of $\zeta \cup \eta$ into $F$, as indicated by the dotted line in Figure~\ref{uglycasefig}.  
\begin{figure}[htb]
  \begin{center}
  \includegraphics[width=1.5in]{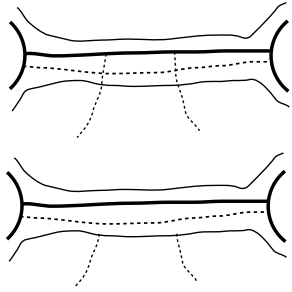}
  \caption{In the case when $\eta$ passes through both $N \times \{\frac{1}{4}\}$ and $N' \times \{\frac{3}{4}\}$, the projection of $\eta \cup \zeta$ may not be embedded, as in the top picture.  We can fix this by pushing $\zeta$ out of $N$, as in the bottom picture.  In the top, the dotted line is $\eta$.  In the bottom, the dotted line is $\sigma$.}
  \label{uglycasefig}
  \end{center}
\end{figure}

Because the arc $\zeta$ is parallel to $\alpha_0$, we can push it into the boundary of $N$ so that $\eta$ is outside $N$ near its endpoints, as follows: The subarc of $\eta$ in $N \times \{\frac{1}{4}\}$ is essential and properly embedded in $F$.  The two arcs in $N' \times \{\frac{3}{4}\}$ are disjoint from $\alpha_0$.  The intersection of the projection of $\eta$ with $N$ consists of arcs $a_1, a_2$ which contain the endpoints of $\eta$.  The union of $a_1$, $a_2$ and the projection of $\zeta$ is a properly embedded arc in $N$ which is parallel to an arc $\sigma$ in $\partial N$.

Let $\sigma'$ be the projection of the complement in $\eta$ of the arcs $a_1$, $a_2$.  Then $\sigma \cup \sigma'$ is a simple closed curve in $F$ which is homotopic to the projection of $\zeta \cup \eta$.  (Note that $\sigma \cup \sigma'$ is not properly embedded because $\sigma'$ intersects $\partial F$.)  Because $\zeta \cup \eta$ bounds the disk $D$, the loop $\sigma \cup \sigma'$ bounds an immersed disk and therefore an embedded disk $E'$.

Consider a short arc $a_3$ properly embedded in $N$ with one endpoint in $\sigma$ and the other endpoint in $\sigma'$.   Such an arc exists because $N$ is connected.  

\begin{Claim}
\label{topclaim2}
If $E'$ contains $a_3$ then $D \cap E$ can be reduced.
\end{Claim}

\begin{proof}
The complement in $E'$ of $a_3$ is a pair of disks.  Sliding the image of $\eta$ across these disks defines an isotopy of $\partial D'$ which brings $\eta$ into $N \times [\frac{1}{4},\frac{3}{4}]$.  From here, $\eta$ can be isotoped disjoint from $E$, reducing the number of intersections.
\end{proof}

\begin{Claim}
\label{topclaim3}
If $E'$ is disjoint from $a_3$ then $D \cap E$ can be reduced.
\end{Claim}

\begin{proof}
The arcs of $\partial D$ adjacent to $\eta$ begin in $E'$.  In order to leave $E'$, they must pass through $\alpha_0$ or enter a component of $\partial F$.  If they pass through $\alpha_0$ in the opposite direction or end at $\partial F$ then the arc cuts off a disk from $E'$.  Sliding the arc across this disk reduces the number of intersections.  If the arc crosses $\alpha$ in the same direction then the next arc ends up in $E'$ again.  In order for the curve to end, there must eventually be an arc which can be isotoped across $\alpha$, reducing $D \cap E$.
\end{proof}

If $D \cap E$ is minimized, Claims~\ref{topclaim1},~\ref{topclaim2} and~\ref{topclaim3} imply that there must be an arc of $\eta$ in $N' \times \{\frac{3}{4}\}$ which is disjoint from $E$, so Claim~\ref{topnclaim} implies that $\alpha_0$ and $\alpha_1$ are disjoint.
\end{proof}

\section{Induction on Almost Bundles}
\label{structsect}

Lemma~\ref{techlem} will now supply the inductive step in the following proof.

\begin{proof}[Proof of Lemma~\ref{spinelem}]
Let $X$ be an almost bundle constructed from a compact, connected, orientable surface $F$, essential arcs $\alpha_0$, $\alpha_1$ and a map $\phi$.  We will induct on the negative Euler characteristic, $-\chi(F)$.   Let $Y$ and $Z$ be defined as in the previous section.  There is an inclusion map $Y \rightarrow X$ which is one-to-one in the interior of $Y$ and sends $Z$ onto the boundary of $X$. 

Because there are essential arcs in $F$, the surface cannot be a disk and $\chi(F) \leq 0$.  If $\chi(F) = 0$ then $F$ is an annulus. Any two properly embedded essential arcs in an annulus are isotopic so we can assume $\alpha_0 = \alpha_1$.  Let $\ell_1$ be an essential loop in the annulus $F \times \{\frac{1}{2}\}$.  Let $\ell_2$ be the loop formed by extending arcs up and down from a point in $\ell_2$ so that the ends of the arcs are identified in $X$.  Then $\ell_1 \cup \ell_2$ is a spine consisting of a vertical loop and a horizontal loop.  This is the base case of the induction argument.

For the inductive step, assume $n = -\chi(F) > 0$ and assume that for any almost bundle $X'$ constructed from a surface $F'$ with $-\chi(F') < n$, $X'$ has a spine of the desired type.  Because $X$ is a handlebody, there is an essential, properly embedded disk $D \subset X$.  We can assume that $D$ is transverse to $F \times \{0\}$ and $D \cap (F \times \{0\}) \subset X$ is a collection of arcs and loops.  

Isotope $D$ so as to minimize the number of components of intersection.  Because $X$ is irreducible and $F$ is incompressible, minimization implies that $D \cap (F \times \{0\})$ is a collection of arcs.  If the intersection is empty, let $D' = D$.  Otherwise, let $D'$ be an outermost disk of the complement $D \setminus (F \times \{0\})$.

The pre-image of $D'$ in $Y$ is properly embedded.  If $\partial D'$ bounds a disk in $\partial Y$ and $D' = D$ then $D$ is boundary parallel in $X$.  If $\partial D'$ is trivial and $D'$ is a proper subdisk of $D$ then isotope $D$ across the disk in $\partial Y$ bounded by $\partial D'$.  This isotopy reduces the number of components of $D \cap (F \times \{0\})$.  Because we assumed $D$ is essential and $D \cap (F \times \{0\})$ is minimal, we know $\partial D'$ is essential in $Y$.  

The complement in $\partial X$ of $Z$ is a regular neighborhood of $\partial X \cap (F \times \{0\})$.  If $D' = D$ then $D'$ is disjoint from $F \times \{0\}$ so $\partial D'$ can be isotoped into $Z$.  Otherwise, after a suitable isotopy, $\partial D' \setminus Z$ is a single arc because $D'$ is an outermost disk.  In either case, $\partial D'$ meets the criteria of Lemma~\ref{techlem}, implying that we can isotope $\alpha_0$ to be disjoint from $\alpha_1$.  

Let $N$ and $N'$ be disjoint regular neighborhoods of $\alpha_0$ and $\alpha_1$, respectively.  Let $N_0 = N \times [0,\frac{3}{4})$ and $N_1 = N' \times (\frac{1}{4},1]$.  Note that we have essentially taken the original sets $N_0$ and $N_1$, and extended them up and down, respectively.  Because $N$ and $N'$ are disjoint outside a regular neighborhood of $\partial F$, the open sets $N_0$ and $N_1$ intersect only within a regular neighborhood of $\partial F \times [0,1]$.  They are regular neighborhoods of $\alpha_0 \times \{0\}$ and $\alpha_1 \times \{1\}$, respectively, so gluing along their complements yields a manifold $X'$ which is homeomorphic to $X$, but has a different fiber structure.

Define $F' = F \times \{\frac{1}{2}\} \setminus N_1$, $\alpha'_0 = \alpha_0 \times \{\frac{1}{2}\} \subset F'$ and $\alpha'_1 = \phi(\alpha_1) \times \{\frac{1}{2}\}$.  We claim that $X'$ is an almost bundle defined by $F'$, $\alpha'_0$ and $\alpha'_1$.  To see this, consider the following construction:

Cut $F \times [0,1]$ along $F \times \{\frac{1}{2}\}$ to produce two pieces and let $\psi : F \times \{\frac{1}{2}\} \rightarrow F \times \{\frac{1}{2}\}$ be the identity map from the ``top'' of $F \times [0,\frac{1}{2}]$ to the ``bottom'' of $F \times [\frac{1}{2},1]$.  If we remove $N_0$ and $N_1$ from $F \times [0,1]$, then cut along $F \times \{\frac{1}{2}\}$, the bottom piece looks like $F' \times [0,\frac{1}{2}]$ with $N' \times (\frac{1}{4},\frac{1}{2}]$ removed.  The top piece is $F' \times [\frac{1}{2},1]$ with $N \times [\frac{1}{2},\frac{3}{4})$ removed.  Now take $\psi$ to be the restriction of the identity to $F'$.

Glue the bottom of $F \times [0,\frac{1}{2}]$ to the top of $F \times [\frac{1}{2},1]$ by the map $\phi$. The result is homeomorphic to $F' \times [0,1]$ with the neighborhoods of the arc $\alpha'_0$ in $F' \times \{0\}$ and the arc $\alpha'_1$ in $F \times \{1\}$ removed.  We can recover $X'$ by gluing $F' \times \{0\}$ to $F' \times \{1\}$ by the map $\psi$.
 
By definition, this will be an almost bundle if $\alpha'_0$ and $\alpha'_1$ are essential in $F'$.  Because there is a homeomorphism of $F'$ sending one to the other, we need only show that $\alpha'_0$ is essential.  This is equivalent to showing that $\alpha_0 \cup \alpha_1$ does not separate a disk from $F$.  There are two ways that $\alpha_0 \cup \alpha_1$ can separate a disk from $F$: either $\alpha_0$ and $\alpha_1$ are parallel or one of the arcs separates from $F$ an annulus in which the other arc is essential.  The second case is impossible because the automorphism would send a non-separating arc to a separating arc (or vice versa.)

If $\alpha_0$ and $\alpha_1$ are parallel then there is a vertical disk in $X$ as in the $\chi(F) = 0$ case.  Because $X$ is a handlebody, the complement of this disk is one or two solid tori, implying that $F \setminus \alpha_0$ is one or two disks.  This is impossible if $\chi(F) < 0$ so we conclude that $\alpha'_0$ is essential and $X'$ is an almost bundle constructed from $F'$.

Because $X'$ is an almost bundle and $\chi(F') > \chi(F)$, we know that there is spine of $X'$ consisting of a vertical loop and a horizontal loop.  Let $l_1 \cup l_2$ be the image in the inclusion map $X' \rightarrow X$ of this spine.  This inclusion map is isotopic to a homeomorphism from  $X$ to $X'$ so $l_1 \cup l_2$ is a spine for $X$.  The incluion map sends fibers to fibers to the image of $l_1$ is horizontal in $X$ and the image of $l_2$ is vertical.
\end{proof}

\section{Proof of the Main Theorem}
\label{mainthmsect}

We have proved Lemma~\ref{spinelem}.  Lemma~\ref{bundlehomeolem} will be proved in the second half of the paper and independently of Lemma~\ref{nicespine1lem}.  For now, assume Lemma~\ref{bundlehomeolem} is true.

\begin{proof}[Proof of Lemma~\ref{nicespine1lem}]
By Lemma~\ref{bundlehomeolem}, we can isotope $H_1$ so that there is a fiber-preserving homeomorphism $h : X \rightarrow H_1$ for some almost bundle $X$.  By Lemma~\ref{spinelem}, there is a spine $K$ for $X$ consisting of a single vertical edge and a single horizontal edge.  Because $h$ preserves the fiber structures, the graph $h(K)$ consists of a vertical edge and a horizontal edge with respect to the bundle structure of $M$.  Because $h$ is a homeomorphism, $h(K)$ is a spine for $H_1$.
\end{proof}

\begin{proof}[Proof of Theorem~\ref{mainthm}]
Let $(\Sigma, H_1, H_2)$ be a genus two Heegaard splitting of a surface bundle $M$ with monodromy $\phi$.  By Lemma~\ref{nicespine1lem}, we can isotope $\Sigma$ so that there is a spine $K$ of $H_1$ consisting of a single horizontal loop $\ell_1$ and a single vertical loop $\ell_2$.  Let $S$ be the horizontal surface containing $\ell_1$.

Let $S$ be the level surface containing $\ell_1$ and $h : S \times [0,1] \rightarrow M$ a map that sends $S \times \{0\}$ and $S \times \{1\}$ to $S$.  By choosing $h$ carefully, we can ensure that for some finite collection $C$ of points in $S$, the preimage in $\psi$ of $\ell_2$ is precisely $C \times [0,1]$.  

Let $O \subset S$ be a regular neighborhood in $S$ of $\ell_1$, $N$ the preimage in $S \times \{0\}$ of $O$ and $N'$ the preimage in $S \times \{1\}$ of $O$.  The monodromy $\phi$ takes $C$ onto itself, so there is a regular neighborhod  $N''$ of $C$ such that $\phi(N'') = N''$.  Define $R \subset M$ to be the image in $h$ of the set $(N \times 	[0,\frac{1}{4})) \cup (N' \times (\frac{3}{4}, 1]) \cup (N'' \times [0,1])$.  The set $R$ is a regular neighborhood of the spine $\ell_1 \cup \ell_2$, so there is an isotopy of $M$ taking $H_1$ onto the closure $H_1^{2g}$ of $R$ and $H_2$ onto the complement $H_2^{2g}$ of $R$.

Let $\alpha_0$ be the preimage in $S \times \{0\}$ of $\ell_1$ and $\alpha_{2g-1}$ the preimage in $S \times \{1\}$ of $\ell_1$.  The complement in $S \times [0,1]$ of $R$ is homeomorphic to $F_{2g} \times [0,1]$ where $F_{2g}$ is the complement in $S$ of a regular neighborhood of $C$.  The pre-image of $H_2^{2g}$ in $S \times [0,1]$ is precisely the complement in $F_{2g} \times [0,1]$ of the neighborhood $N \times [0,\frac{1}{4})$ of $(\alpha_0 \times \{0\})$ and the neighborhood $N' \times (\frac{3}{4},1]$ of $(\alpha_{2g-1} \times \{1\})$.  Therefore the bundle structure of $M$ makes $H_2^{2g}$ an almost bundle.

By Lemma~\ref{techlem}, there is a homeomorphism $\psi_{2g}: F_{2g} \rightarrow F_{2g}$, isotopic to the identity, such that $\psi_{2g}(\alpha_{2g-1})$ is disjoint from $\alpha_0$.  We can choose $\psi_{2g}$ to restrict to the identity on the boundary of $F_{2g}$, although the isotopy may not fix $\partial F_{2g}$.  This $\psi_{2g}$ extends to a map from $S$ to $S$  which is isotopic to a sequence of Dehn twists along $\partial F_{2g}$.  The boundary of $F_{2g}$ consists of a number of trivial loops in $S$, so the extension of $\psi_{2g}$ is isotopic to the identity on $S$.  If we replace $\alpha_{2g-1}$ with $\psi_{2g}(\alpha_{2g-1})$ then $\alpha_{2g-1}$ will be disjoint from $\alpha_0$ and $\psi_{2g}^{-1} \circ \phi$ will send $\alpha_{2g-1}$ onto $\alpha_0$.

Let $M^{2g-1}$ be the surface bundle defined by the monodromy $\psi_{2g}^{-1} \circ \phi$.  This manifold is homeomorphic to $M$ and there is a Heegaard splitting induced from $M$, which we will also denote $(\Sigma^{2g}, H_1^{2g}, H_2^{2g})$, and a map $S \times [0,1] \rightarrow M^{2g-1}$ such that the horizontal arc of the spine of $H^{2g}_1$ $\alpha_0 \times \{0\}$ and $\alpha_{2g-1} \times \{1\}$ whose projections into $S$ are disjoint.  Let $D_{2g-1}$ be the image in $M^1$ of  $\alpha_{2g-1} \times [0,1]$ and let $H^{2g-1}_1$ be the closure of a regular neighborhood $N$ of $K \cup D_{2g-1}$.  Let $H_2^{2g-1}$ be the complement of $N$.

The disk $\alpha_{2g-1} \times [0,1]$ is disjoint from $\alpha_0 \times \{0\}$ and the monodromy map sends $\alpha_0$ onto $\alpha_{2g-1}$ so the image of $\alpha_{2g-1}$ in the monodromy is an arc $\alpha_{2g-2}$ which is disjoint from $\alpha_{2g-1}$.  Let $F_{2g-1}$ be the complement in $F_{2g}$ of an open neighborhood of $\alpha_{2g-1}$.  The preimage of $H_2^{2g-1}$ is the complement in $F_{2g-1} \times [0,1] (\subset S \times [0,1])$ of a neighborhood of $\alpha_0 \times \{0\}$ and an open neighborhood of $\alpha_{2g-2} \times \{1\}$.  So, $H_2^{2g-1}$ is again an almost bundle.

As before, Lemma~\ref{techlem} implies that there is a homeomorphism $\psi_{2g-1} : F_{2g-1} \rightarrow F_{2g-1}$ which is isotopic to the identity on $F_{2g-1}$ and extends to a map from $S$ to $S$.  The map $\psi_{2g-1} : S \rightarrow S$ is isotopic to a number of dehn twists along $\partial F_{2g-1}$.  In this case, $\partial F_{2g-1}$ consists of a number of trivial loops plus the boundary of a regular neighborhood of $\alpha_{2g-1}$.  

Replace $\alpha_2$ with $\psi_{2g-2}(\alpha_2)$, so that the arcs $\alpha_0$, $\alpha_1$, $\alpha_2$ will be pairwise disjoint and $\psi_{2g-1}^{-1} \circ \psi^{-1}_{2g} \cup \phi$ sends $\alpha_{2g-2}$ onto $\alpha_{2g-1}$ and $\alpha_{2g-1}$ onto $\alpha_0$.  Let $(\Sigma^{2g-2},H_1^{2g-2}, H_2^{2g-2})$ be the Heegaard splitting formed from $(\Sigma^{2g-1},H_1^{2g-1}, H_2^{2g-1})$ by extending $H_1^{2g-1}$ along the disk $D_{2g-2} = \alpha_{2g-2} \times [0,1]$.

We again construct a manifold $M^{2g-2}$, with monodromy $\psi^{-1}_{2g-1} \circ \psi^{-1}_{2g} \circ \phi$ and induced Heegaard splitting $(\Sigma^{2g-2}, H_1^{2g-2}, H_2^{2g-2})$.  We can continue the process of forming a new manifold $M_{2g-i}$ for each $i < 2g$ and then extending the induced spine $K \cup D_{2g-1} \cup \dots \cup D_{2g-i+1}$ by a disk $D_{2g-1} = \alpha_{2g-i} \times [0,1]$.  The process will terminate when the arcs $\alpha_0$ and $\alpha_{2g-i}$ are isotopic in $F_{2g-i}$.  Because $H_2^{2g-i}$ is a handlebody, this can only occur when $F_{2g-i} \setminus \alpha_0$ is a disk, so $F_{2g-i} = F_1$ is an annulus.

Construct one final manifold $M^0$ with monodromy $\psi = \psi_{1}^{-1} \circ \dots \circ \psi_{2g} \circ \phi$.  The image in $M^0$ of each disk $D_i$ intersects $S$ in the arcs $\alpha_{i-1}$ and $\alpha_i$.  Because $D_i$ is isotopic to $\alpha_{i-1} \times [0,1]$ in $S \times [0,1]$, the map $\phi$ must send $\alpha_{i-1}$ to $\alpha_i$.   Thus $\psi$ is a cyclic rotation of the one-vertex spine $\bigcup \alpha_i$ of $S$.  

By composing both sides with each $\psi_i$, we find that $\psi_{2g} \circ \phi = \psi_{1} \circ \dots \circ \psi_{2g-1} \circ \psi_0$.  Recall that for $i, j > 0$, $\psi_i$ and $\psi_j$ are isotopic to Dehn twists along disjoint loops so $\phi_i$ and $\psi_j$ commute up to isotopy.  Also, $\phi_{2g}$ is isotopic to the identity on $S$ so we find that $\phi$ is isotopic to $\psi_{2g-1} \circ \dots \circ \psi_{0}$.
\end{proof}

\section{General Position}
\label{positionsect}

Our goal now is to isotope $\Sigma$ so that the fiber structure of the fiber bundle defines an almost bundle structures on the handlebodies of the Heegaard splitting.  We will progressively refine the position of the Heegaard surface, motivated by the following definitions.

\begin{Def}
A surface $\Sigma$ in a surface bundle $M$ with bundle map $\pi : M \rightarrow S^1$ is \textit{essentially embedded} if the restriction $\pi|_{\Sigma}$ is a (circle-valued) Morse function and each loop in each level set is essential in both $\Sigma$ and in the level surface $\pi^{-1}(x)$.
\end{Def}

Note that the definition refers to the loops of the level sets, not the components of the level sets.  The components of level sets that contain critical points will not be loops.  Also note that an essentially embedded surface may be compressible.  However, every essential surface is isotopic to an essentially embedded surface by a theorem in~\cite{thr:norm}.

\begin{Def}
A surface $\Sigma$ in a surface bundle $M$ with bundle map $\pi : M \rightarrow S^1$ is \textit{fair} if the restriction $\pi|_{\Sigma}$ is a (circle-valued) Morse function and each loop in each level set is either essential in both $\Sigma$ and in the appropriate level surface of $\pi$ or trivial in both surfaces.
\end{Def}

We will say that a loop of intersection is \textit{fair} if it is either essential in both surfaces or trivial in both.  Thus a surface is fair if each transverse intersection with a level surface is fair.

Every essentially embedded surface is fair.  Also, every incompressible surface which is in Morse position is fair because the level surfaces of the bundle are incompressible.  A loop which is essential in one surface and trivial in the other defines a compression for one of the surfaces.  

We would like to have the Heegaard surface be essentially embedded.  Unfortunately, in general, we need a slightly weaker condition.  This is reminiscent of the fact that a strongly irreducible Heegaard splitting can be made almost normal in a triangulation, but not necessarily normal.

\begin{Def}
A surface $\Sigma$ in a surface bundle $M$ with bundle map $\pi : M \rightarrow S^1$ is \textit{almost fair} if $\Sigma$ is fair or $\Sigma$ is the result of attaching a horizontal handle to a fair surface.
\end{Def}

\begin{figure}[htb]
  \begin{center}
  \includegraphics[width=3.5in]{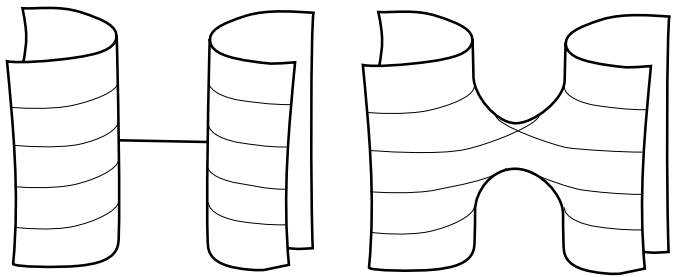}
  \caption{Attaching a horizontal handle to a surface.}
  \label{addhandlefig}
  \end{center}
\end{figure}

By attaching a horizontal handle, we mean the following:  Let $\Sigma$ be a surface in a surface bundle $M$ and let $\alpha$ be an arc in a level surface whose endpoints are in $\Sigma$ and whose interior is disjoint from $\Sigma$.  Let $N$ be a regular neighborhood of $\alpha$ in $M$.  Then $\partial N \setminus \Sigma$ consists of two disks and an annulus $A$.  Let $\Sigma'$ be the result of removing the interior of $N$ from $\Sigma$ and attaching $A$.  The construction is shown in Figure~\ref{addhandlefig}

In the following sections, we will show that every strongly irreducible Heegaard splitting of a surface bundle is isotopic to an almost essential surface and every genus two Heegaard splitting is isotopic to an essentially embedded surface.  The first few sections of the proof follow arguments of Bachman and Schleimer~\cite{bsc:bndls}.

\section{Sweep-outs}
\label{sweepsect}

We will first show that every strongly irreducible Heegaard splitting is isotopic to an almost fair surface.  To do this, we will use the double sweep-out method introduced by Rubinstein and Scharlemann~\cite{rub:compar}.

Let $M$ be a 3-manifold and $(\Sigma, H_1, H_2)$ a Heegaard splitting for $M$.   Let $K_1$ be a spine of $H_1$ and let $K_2$ be a spine of $H_2$.  The complement in $M$ of $K_1 \cup K_2$ is homeomorphic to $\Sigma \times (-1,1)$.  

The projection of this homeomorphism onto the $(-1,1)$ factor extends to a map $f : M \rightarrow [-1,1]$ such that $f^{-1}(-1) = K_1$, $f^{-1}(1) = K_2$ and for every $x \in (-1,1)$, $f^{-1}(x)$ is a surface isotopic to $\Sigma$.  The map $f$ is a called \textit{sweep-out} of $(\Sigma, H_1, H_2)$ and we can choose $f$ to be a smooth function.

Let $\pi : M \rightarrow S^1$ be a map that defines a surface bundle structure on $M$.  The \textit{discriminant set} or \textit{Jacobi set} of $(f, \pi)$ is the set $J = \{x \in M : \nabla f(x) = \lambda \nabla \pi(x)$ for some $\lambda \in \mathbf{R} \}$.  Equivalently, $J$ is the set of points in $M$ where a level surface of the sweep-out is tangent to a level surface of the bundle structure or where $\nabla f = 0$.   (Because $f$ is smooth, the set of points where $\nabla f$ is 0 is precisely the two spines of the handlebodies.)

A slight modification of the proof used by Kobayashi~\cite{Kob:disc} for a pair of sweep-outs implies that after an arbitrarily small isotopy, the spines can be made transverse to the level surfaces of $\pi$ at all but finitely many points.  A further isotopy ensures that $J$ is a one-dimensional set consisting of the spines of $H_1$ and $H_2$ (where $\nabla f = 0$) and a number of edges and loops in their complement.  Each edge in the complement has both endpoints in the spines.  Let $G = (f \times \pi)$ be the map $M \rightarrow [-1,1] \times S^1$ defined by $G(p) = (f(p),g(p))$.  The \textit{Rubinstein-Scharlemann graphic} is the image $R = G(J)$ in $[-1,1] \times S^1$.  

We can also interpret the Rubinstein-Scharlemann graphic as follows:  For each level surface $f_x = f^{-1}(x)$ of the sweep-out, the restriction $\pi_x = \pi|_{f_x}$ is a map from $f_x$ to $S^1$.  The set $J \cap f_x$ is precisely the set of critical points of $\pi_x = 0$ which coincide with tangencies of the level surfaces of $f$ and $\pi$.  When $f$ and $\pi$ are in general position, $\pi_x$ is a circle valued Morse function for all but finitely many values of $x$.  The image $G(f_x)$ is the circle $\{x\} \times S^1$ and this circle intersects the graphic at the critical values of $\pi_x$.  As the circle sweeps through a small interval around $x$, the images of the critical values form edges in the graphic.  

Vertices in the graphic result from the finitely many values of $x$ for which $\pi_x$ is not a proper Morse function.  There are two possible ways that $\pi_x$ can fail to be a proper Morse function: There may be two critical points at the same level or there may be a single degenerate critical point.  We will call such functions \textit{near Morse} functions.

We will think of the family $\{\pi_x : x \in (-1,1)\}$ as a path in the space of functions on $\Sigma$.   When this path passes through a patch where two critical points pass each other, it must pass through a function where two critical points are at the same level, creating a valence-four vertex in the graphic.  General position implies that in any horizontal arc $[0,1] \times \{y\}$ or vertical circle $\{x\} \times S^1$, there is at most one crossing.  If the path passes through a patch where two critical points cancel or where two critical points are created (uncanceled?), there is a function with a degenerate critical point, creating a cusp in the graphic.

A \textit{region} of the graphic is a component of the complement $([-1,1] \times S^1) \setminus R$.  For a region $C$ of the graphic, let $(x,y) \in C$ and  let $L_C = G^{-1}(x,y)$.  This set is a collection of loops, the intersection of the surfaces $f^{-1}(x)$ and $\pi^{-1}(y)$.  Given a second point, $(x',y')$ in the same region as $(x,y)$, a piecewise vertical and horizontal path from $(x,y)$ to $(x',y')$ defines an (ambient) isotopy of $M$ taking $f^{-1}(x)$ to $f^{-1}(x')$ and $g^{-1}(y)$ to $g^{-1}(y')$.  Thus up to isotopy, the set $L_C$ depends only on the component $C$, not on the choice of $(x,y)$.

Label each region of the graphic as follows:  If there is a loop in $L_C$ which is essential in $\Sigma$ and bounds a disk in $f^{-1}([-1,x])$, then label the region with a 1.  If there is a loop in $L_C$ which is essential in $\Sigma$ and bounds a disk in $f^{-1}([x,1])$, label the region with a 2.  Otherwise, leave the region unlabeled.  

\begin{Def}
A Heegaard splitting $(\Sigma, H_1, H_2)$ is \textit{strongly irreducible} if for every pair of properly embedded disks $D_1 \subset H_1$ and $D_2 \subset H_2$, we have $\partial D_1 \cap \partial D_2 \subset \Sigma$ is not empty.  If $\Sigma$ is not strongly irreducible then $\Sigma$ is \textit{weakly reducible}.
\end{Def}

If $L_C$ contains a loop which is trivial in $F$, but essential in $\Sigma$ then by Scharlemann's no nesting Lemma~\cite{sch:detect}, it contains a loop bounding a disk in one of the handlebodies or $\Sigma$ is weakly reducible.  There cannot be a loop which is essential in $F$ but inessential in $\Sigma$ because such a loop would suggest a compressing disk for $F$.  Thus if $\Sigma$ is strongly irreducible then in an unlabeled region, each curve of intersection is either essential in both $F$ and $\Sigma$ or trivial in both.

\section{Analyzing the Graphic}

Let $M$, $\pi$, and  $(\Sigma, H_1, H_2)$ be defined as in the previous section and assume that the level surfaces of $\pi$ have genus two or greater.

\begin{Lem}
\label{splitlem}
If $\Sigma$ is a strongly irreducible Heegaard splitting then $\Sigma$ is isotopic to an almost fair surface.
\end{Lem}

\begin{proof}
Let $f$ be a sweep-out of $(\Sigma, H_1, H_2)$ and assume $f, \pi$ are in general position. Let $(x,y)$ be a point in a region $C$ of the graphic.  The level surface $f_x$ of $f$ is the preimage in $G$ of a vertical circle in $\{x\} \times S^1$. The loops in $L_C = f_x^{-1}(y)$ are level curves of the Morse function $\pi_x$ on the surface $f_x$ so loops that do not coincide are disjoint.  Thus if the vertical circle passes through a region labeled with a 1 and a region labeled with a 2 (or a region with both labels,) then the Heegaard splitting is weakly reducible.  Because we assumed $\Sigma$ is strongly irreducible, no region can have both labels and if a vertical circle passes through two labeled regions, they must have the same label.

For a generic level surface $S = \pi^{-1}(y)$, the restriction of $f$ to $S$ is a Morse function so there is a value $\varepsilon$ such that $G^{-1}([-1,-1+\varepsilon],y)$ consists of a number of disks containing only index-zero critical points.  Thus all the loops of intersection between $S$ and $f_{-1+\varepsilon}$ bound disks in $f^{-1}([-1, -1+\varepsilon])$.  Each level surface of $\pi$ must intersect each Heegaard surface because handlebodies do not contain incompressible surfaces.  Thus any region adjacent to the circle $\{-1\} \times S^1$ is labeled with a 1.  A similar argument implies that any region adjacent to $\{1\} \times S^1$ is labeled with a 2.

Because no circle can pass through regions with different labels, there must be a vertical circle $\{x\} \times S^1$ which does not pass through any labeled regions.  If this circle is disjoint from all the vertices of the graphic then the corresponding function $\pi_x$ on $f_x$ is Morse and all the loops of intersection between $\Sigma$ and level surfaces of $\pi$ are essential in both or trivial in both.  In other words, the surface $f_x$ is fair and the proof is complete.

Assume that $\{x\} \times S^1$ intersects one or more vertices of $R$.  Each vertex that $\{x\} \times S^1$ intersects corresponds to a degenerate critical point in $\pi_x$ or a level with two critical points.  Because $\pi_x$ is Morse or almost Morse, the vertical circle can pass through at most one vertex.  If the region to the left or the right of the vertex is unlabeled, then a vertical circle slightly to the left or the right will pass through all unlabeled regions and no vertices, and again the proof would be complete.  Thus we can assume that the region to the left is labeled with a 1, the region to the right is labeled with a 2 and the regions above and below are unlabeled. 

Assume for contradiction this is the case. Let $B$ be the region labeled with a 1, let $D$ be the region labeled with a 2 and let $A$ and $C$ be the remaining regions.  We will think of $A$ as being above the vertex, $B$ to the left, $C$ below it and $D$ to the right, as in Figure~\ref{crossingfig}.  The label of a region does not change across an edge consisting of central singularities so the two edges involved in the crossing must both be saddle singularities.
\begin{figure}[htb]
  \begin{center}
  \includegraphics[width=1.5in]{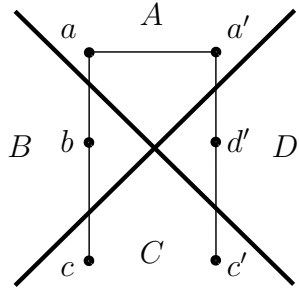}
  \put(-60,100){$A$}
  \put(-110,50){$B$}
  \put(-60,10){$C$}
  \put(-10,50){$D$}
  \put(-90,95){$a$}
  \put(-27,95){$a'$}
  \put(-90,52){$b$}
  \put(-27,52){$d'$}
  \put(-90,5){$c$}
  \put(-27,5){$c'$}
  \caption{The setup for analyzing a vertex where the label changes.}
  \label{crossingfig}
  \end{center}
\end{figure}

Let $\alpha \subset [-1,1] \times S^1$ be a short vertical arc which passes from region $A$ to region $B$, then to region $C$.  Let $a$ be the endpoint of $\alpha$ in $A$, $b$ a point of $\alpha$ in $B$ and $c$ the endpoint in $C$, as in Figure~\ref{crossingfig}.  The preimage in $G$ of the arc $[a,b]$ is a subsurface of $f_x$.  Exactly one component of this surface contains a critical point, and this critical point is a saddle.  Let $X_0$ be this subsurface.  Similarly, let $X_1$ be the component of the preimage of $[b,c]$ containing a saddle.

\begin{figure}[htb]
  \begin{center}
  \includegraphics[width=3.5in]{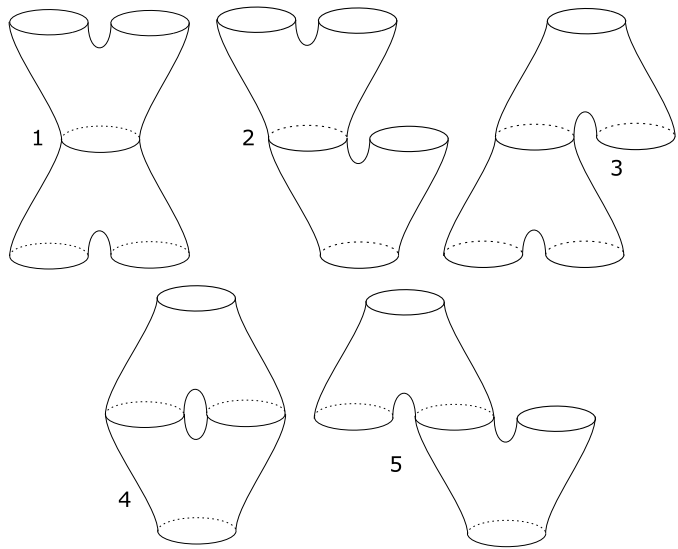}
  \caption{Possible configurations for the two pairs of pants, $X_0$, $X_1$, in the preimage of the arc $ac$.}
  \label{twopantsfig}
  \end{center}
\end{figure}
If $X_0$ and $X_1$ do not share a boundary component, then the isotopy classes of the level loops will not change when the critical levels pass each other.  Because the label changes at the vertex, the subsurfaces $X_0$ and $X_1$ must share at least one boundary component.  The possible configurations of $X_0 \cup X_1$ are shown in Figure~\ref{twopantsfig}.   We will rule out the first four configurations and show that the last configuration allows us to make $\Sigma$ into an almost fair surface.

Let $\alpha'$ be a vertical arc to the right of the vertex, starting in region $A$, passing through $D$ and into $C$.  Let $a'$, $d'$ and $c'$ be points on $\alpha'$ and assume $a$ and $a'$ are on the same level in $S^1$, as are $c$ and $c'$.  As with $\alpha$, the preimages of the arcs $[a',d']$ and $[d',c']$ contain pairs of pants $X'_0$ and $X'_1$, respectively.  

When $f_t$ passes through the crossing, the isotopy classes in $f_t$ of the boundary components of $X_0 \cup X_1$ do not change. The loop (or loops) $X_0 \cap X_1$, however, is replaced by the not necessarily isotopic $X'_0 \cap X'_1$.  The set $\partial X_0 \cup X_1$ is isopic to $\partial X'_0 \cup X'_1$, so the configuration of $X'_0 \cup X'_1$ is the same as the configuration of $X_0 \cup X_1$.  This also implies that the boundary components of $X_0 \cup X_1$ are all fair loops.  Because regions $A$ and $B$ are labeled, $\partial X_0 \cap \partial X_1$ must contain a loop which is essential in $\Sigma$ but trivial in $S$, as must $\partial X'_0 \cap \partial X'_1$.

Let $\beta$ be a horizontal arc from $a$ to $a'$.  For each boundary component at the top of $X_0$, the preimage of $\beta$ contains an annulus from this loop to the corresponding boundary component of $X'_0$.  Likewise, the preimage of an arc $\beta'$ from $c$ to $c'$ contains annuli connecting the bottom boundary components of $X_1$ and $X'_1$.

If $X_0 \cup X_1$ is in configuration 1 or 2 then the union of $X_0$, $X'_0$ and the appropriate annuli from the preimage of $\beta$ form a twice punctured torus.  The punctures are the loops $X_0 \cap X_1$ and $X'_0 \cap X'_1$, which bound disks in $S$.

There is an isotopy in $[-1,1] \times S^1$, transverse to the graphic, taking the arc $[a,b]$ to a horizontal arc which intersects $R$ in a single point.  Kobayashi~\cite{Kob:disc} showed that such an isotopy in the graphic induces an isotopy in $M$.  This isotopy takes $X_0$ to a pair of pants in a level surface of $\pi$.  Because a loop in $X_0 \cap X_1$ bounds a disk in $F$, so does the corresponding boundary component of the level pair of pants. 

Similarly, $X'_0$ is isotopic to a horizontal pair of pants with a trivial boundary component.  The preimage of the two horizontal arcs and the arc $[a,a']$ form a twice-punctured torus in a level surface such that the boundary components are trivial in $S$. This implies that the level surface is a torus.  Because we assumed $S$ is not a torus, this rules out configurations 1 and 2.  A similar argument for $X_1$ and $X'_1$ rules out configurations 1 and 3, so we are left with the last two configurations.

In configuration 4, as $f_x$ passes through the vertex, the Morse function changes by replacing a level loop, $X_0 \cap X_1$ with a new level loop $X'_0 \cap X'_1$ so that each component of $X_0 \cap X_1$ intersects each loop of $X'_0 \cap X'_1$ in a single point.  Because at least one loop of each intersection bounds a disk and these disks are on opposite sides of $f_t$, this implies that $(\Sigma, H_1, H_2)$ is stabilized and rules out configuration 4.

For configuration 5, consider what happens to $X_0 \cup X_1$ if we slide the arc $\alpha$ onto the vertex in the graphic.  The two saddles come together to form three loops connected at two points, as in Figure~\ref{hdiskfig}.  The middle loop bounds a disk $D$ in $F$ and the remaining loops are all essential in $F$.  Because all the other level sets in $\Sigma$ are fair, if $\Sigma$ intersects the interior of $D$, then the loop of intersection is trivial in $\Sigma$ and $\Sigma$ can be isotoped disjoint from the interior of $D$.
\begin{figure}[htb]
  \begin{center}
  \includegraphics[width=1.5in]{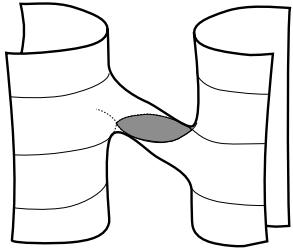}
  \caption{At the crossing, the level surface $f_t$ has two saddle singularities at the same the singular component of the critical level bounds a disk in one of the bundle surfaces, implying that $f_t$ is almost fair.}
  \label{hdiskfig}
  \end{center}
\end{figure}

A slight isotopy of $D$ makes it transverse to the level surfaces of $\pi$.  If we  collapse $\Sigma$ along this transverse disk, the resulting surface, $\Sigma'$, is fair and we can recover $\Sigma$ from $\Sigma'$ by attaching an annulus that undoes the collapse along $D$.  Because $D$ is transverse to the fibration, $\Sigma$ is isotopic to the result of adding a handle to $\Sigma'$ along a horizontal arc.  Thus if $\Sigma$ is not fair, it is almost fair.
\end{proof}

\section{Almost Fair Surfaces}

We will now restrict our attention to genus two Heegaard surfaces to show that an almost fair genus two Heegaard surface must in fact be fair.  We will begin by examining the case of a fair torus.

\begin{Lem}
\label{fairtoruslem}
If $T$ is a fair torus in a surface bundle $M$ then $T$ is isotopic to either an essentially embedded torus in $M$ or the boundary of a regular neighborhood of an essential horizontal loop.
\end{Lem}

\begin{proof}
Choose $x_1,\dots,x_n \in S^1$ so that the disjoint union of surfaces $F = \pi^{-1}(\{x_1,\dots,x_n\})$ cuts $T$ into planar pieces. This can be done, for example by choosing an $x_i$ between every pair of consecutive critical points in $\pi|_\Sigma$.   Assume $T$ has been isotoped so as to minimize $F \cap T$ while preserving the properties that $T$ is a fair surface and $F \setminus T$ consists of planar surfaces.  

The complement $M \setminus F$ is homeomorphic to $F \times (0,1)$, a disjoint union of copies of $S \times (0,1)$ and we will fix a homeomorphism.  Let $\Sigma'$ be the closure in $F \times [0,1]$ of the image of $\Sigma \setminus F$.  Then every component of $\Sigma'$ is a closed planar surface properly embedded in $F \times [0,1]$.

Because $F$ is a closed surface, every properly embedded disk in $F \times [0,1]$ is boundary parallel.  Thus if some component of $\Sigma'$ is a disk, the component can be pushed into $F$.  The corresponding isotopy of $\Sigma$ reduces the number of intersections with $F$ while preserving the two required properties.  Because the number of intersections is minimal, no component of $\Sigma'$ can be a disk.  Any planar surface in a torus whose complement does not contain a disk must be an annulus with essential boundary loops.  Thus each component of $\Sigma'$ must be an annulus  and every loop in $F \cap \Sigma$ must be essential in $\Sigma$.  Because $\Sigma$ is fair, this implies every loop is essential in $F$ as well.

Every essential annulus in $F \times [0,1]$ is isotopic to a vertical annulus, i.e. an annulus of the form $\ell \times [0,1]$ for some simple closed curve $\ell$ in $F$.  Any annulus in $\Sigma'$ which is not vertical is either compressible or boundary compressible.  If an annulus component $A$ of $\Sigma'$ is compressible then the result of compressing $A$ is a pair of disks.  The boundaries of these disks are trivial in $F$, so the boundary components of $A$ are trivial in $F$.  This is impossible because we showed that every loop of $F \cap A$ must be essential in $F$.

An annulus in $\Sigma$ which is not vertical or compressible must be boundary compressible.  Let $A \subset \Sigma'$ be a boundary compressible annulus in $F \times [0,1]$ and let $D$ be the disk that results from boundary compressing $A$.  Let $D' \subset F$ be the disk such that $\partial D' = \partial D$.  We can recover $A$ from $D$ by attaching a band to $D$ along an arc $\alpha \subset F$.  If $\alpha$ were to sit inside $D'$ then the boundary components of $A$ would be trivial loops.  Thus $\alpha$ must be disjoint from the interior of $D'$.  The resulting annulus is therefore boundary parallel in $F \times [0,1]$.

A boundary parallel annulus in $\Sigma$ can be pushed across $F$ while keeping $\Sigma$ a fair surface.  However, after the isotopy $\Sigma'$ may no longer consist of planar pieces.  If the boundary components of this annulus are shared by different annuli in $\Sigma$ then the resulting piece will be planar.  Thus the only situation in which a boundary parallel annulus  can not be pushed across $F$ is when $\Sigma'$ consists of exactly two boundary parallel annuli.

If $\Sigma \cap F$ is minimized and $\Sigma'$ consists of two boundary parallel annuli then $\Sigma$ is the boundary of a regular neighborhood of an essential loop in $F$.  Otherwise, $\Sigma'$ consists entirely of essential vertical annuli and $\Sigma$ is essentially embedded.
\end{proof}

Although we will not apply Lemma~\ref{fairtoruslem} directly in the next Lemma, the method used to prove Lemma~\ref{fairtoruslem} will be instructive in understanding the proof.

\begin{Lem}
\label{almstfairlem}
If $\Sigma$ is an almost fair, genus two Heegaard surface in a surface bundle $M$ then $\Sigma$ is fair.
\end{Lem}

\begin{proof}
Let $\Sigma$ be an almost fair, genus two Heegaard surface.  If $\Sigma$ is fair, then the proof is complete.  Otherwise, $\Sigma$ is the result of attaching a handle to a fair surface $\Sigma'$ along a horizontal arc $\alpha$.  Let $S$ be the level surface that contains $\alpha$.  Because $\Sigma$ is not fair, $S \setminus (\Sigma' \cup \alpha)$ must contain a disk component $D$ such that the boundary of $D$ contains $\alpha$.  Because $\Sigma$ bounds a genus two handlebody on either side, $\Sigma'$ must bound a solid torus or two solid tori on one side.

The disk $D$ may bound $\alpha$ on only one side or on both sides.  If $D$ sits on one side of $\alpha$ then a meridian for $\alpha$ intersects $D$ in a single point.  This implies that the Heegaard surface $\Sigma$ is stabilized.  Because $\Sigma$ is a genus two surface and $M$ is a surface bundle (and not $S^1 \times S^2$), this is impossible.  Thus $D$ must be on both sides of $\alpha$, so $\alpha$ sits in an annular component of $S \setminus \Sigma'$.  In addition, the endpoints of $\alpha$ must sit in essential loops of $\Sigma \cap S$: If $\alpha$ is attached to a loop bounding a disk $D'$ in a level surface, then $D'$ is essential in the complement of $\Sigma$ and can be made disjoint from $D$, implying that $\Sigma$ is weakly reducible.  A weakly reducible, genus two Heegaard splittings is reducible, so this is impossible.

Let $F$ be a collection of level surfaces which are disjoint from $\alpha$ and cut $\Sigma'$ into planar pieces.  Assume $\Sigma' \cup \alpha$ has been isotoped to minimize $F \cap \Sigma'$ while keeping $\Sigma'$ fair, $\alpha$ level and $\Sigma' \setminus F$ planar.  As in the previous argument, any disk in $\Sigma' \setminus F$ is boundary parallel.  Any level set of $\pi$ contained in a disk is trivial in $\Sigma$.  Because the endpoints of $\alpha$ sit in essential loops, such a disk cannot contain an endpoint of $\alpha$, so minimization implies that there are no disks.

Any annulus in $\Sigma' \setminus F$ is either vertical or boundary parallel.  The only situations in which a boundary parallel annulus cannot be removed is when it contains an endpoint of $\alpha$ and when $\Sigma' \setminus F$ consists of two boundary parallel annuli.  Thus after minimization, $\Sigma' \setminus F$ consists of zero or two horizontal annuli and zero or more vertical annuli.

If $\Sigma' \setminus F$ contains no horizontal annuli, then $\Sigma'$ consists of one or two essentially embedded annuli.  The complement of $\Sigma'$ consists of two or three punctured surface bundles.  If the boundary of a punctured surface bundle is compressible, then the surface is a disk.  Because $\Sigma'$ is essentially embedded, no component of the complement can be a disk boundle.  Thus each component has incompressible boundary, contradicting the assumption that some component of the complement is a solid torus.

Assume $\Sigma' \setminus F$ consists of two horizontal annuli and zero or more vertical annuli.  Then $\Sigma'$ bounds a solid torus $H$.  The intersection $H \cap F$ is a collection of parallel annuli which cut $H$ into a collection of solid tori.  Each component of $H \cap F$ runs along the longitudes of the solid tori, so each curve of intersection between $\Sigma'$ and a level surface is a longitude of $H$.  

We can recover $\Sigma$ from $\Sigma'$ by attaching a handle along $\alpha$.  Because $\alpha$ sits in an annular component of $S \setminus \Sigma'$, we can then compress $\Sigma$ along the complement of $\alpha$, forming a surface $\Sigma''$.  This construction is equivalent to attaching a neighborhood of an annulus to the solid torus along a pair of longitudes.  The resulting pair of tori bounds a component homeomorphic to $T^2 \times I$ and containing $\Sigma$.  Because $\Sigma''$ is a compression of $\Sigma$, the other components must be solid tori.  This implies that $M$ is a lens space, contradicting the assumption that $M$ is a surface bundle (and not $S^1 \times S^2$).  The contradiction implies that $\Sigma$ must be essentially embedded.
\end{proof}

\section{Fair Surfaces}
\label{fairsect}

Let $M$ be a surface bundle, $\pi$ a map defining the bundle structure and $\Sigma$ a genus two Heegaard surface in $M$ such that the restriction of $\pi$ to $\Sigma$ is Morse.  Let $c \in S^1$ be a critical level corresponding to a saddle singularity $v$ of $\pi|_\Sigma$.  The pre-image of a small neighborhood of $c$ contains a thrice-punctured sphere component in $\Sigma$ which contains $v$.  If all three boundary components are essential in $\Sigma$ then $v$ is called an \textit{essential saddle}.  If the boundary components are also essential in the appropriate level sets of $\pi$ then we call the thrice-punctured sphere an \textit{essential pair of pants}.  

Every essential pair of pants can be constructed (up to isotopy) by the following construction: Let $\ell \times [0,1]$ be an essential annulus in $F \times [0,1]$.  (Here, $\ell$ is a simple closed curve in $F$.)  Let $\alpha$ be an arc in $F \times \{0\}$ (or $F \times \{1\}$) with endpoints in $\ell \times \{0\}$ and $N$ a regular neighborhood of $(\ell \times [0,1]) \cup \alpha$.  One component of $\partial N$ is an essential pair of pants.  If there is a second component of $\partial N$, this is an annulus isotopic to $\ell \times [0,1]$.

In other words, an essential pair of pants is the result of attaching a horizontal band to a vertical annulus.  We could have constructed the same essential pair of pants by beginning with two essential annuli and attaching them by a horizontal band in $F \times \{1\}$ ($F \times \{0\}$, respectively).

A genus two surface contains exactly two essential saddles.  If $\Sigma$ is fair then these saddles define essential pairs of pants.  The surface $\Sigma$ is essentially embedded if and only if some finite set of level surfaces of $\pi$ cut $\Sigma$ into essential pairs of pants and vertical annuli.

\begin{Lem}
\label{fairlem}
If $\Sigma$ is a fair genus two Heegaard surface in a surface bundle $M$ then $\Sigma$ is isotopic to an essentially embedded surface.
\end{Lem}

\begin{proof}
Let $x_1, \dots, x_n \in S^1$ and let $F = \pi^{-1}(\{x_1,\dots,x_n\})$.  Choose these points so that the two essential saddles in $\Sigma$ are in essential pairs of pants in different components of $M \setminus F$ and every component of $\Sigma \setminus F$ is a disk, an annulus or a pair of pants.  Isotope $\Sigma$ so that $F \cap \Sigma$ is minimized while these conditions are preserved.  As before, we will define $\Sigma'$ to be the closure of $\Sigma \setminus F$ in $F \times [0,1]$.

As in the last two proofs, any disk component of $\Sigma'$ is boundary parallel in $F \times [0,1]$ so by minimizing $\Sigma'$, we eliminate all disks.  Similarly, any annular component $A$ of $\Sigma'$ is either vertical and essential or has the property that $A$ is horizontal but pushing $A$ through $F$ would create a non-planar component of $\Sigma \setminus F$.

If all the annular components of $\Sigma'$ are vertical then $\Sigma$ is essentially embedded and we're done.  Assume this is not the case.  We have seen that an annular component that is not vertical must be horizontal.  If a horizontal annulus $A$ can not be removed then it shares a boundary component with one of the essential saddles.  

Assume that $\Sigma \setminus F$ contains at least one horizontal annulus and every horizontal annulus in $\Sigma \setminus F$ is adjacent to a pair of pants.  We will show that in this case, $\Sigma$ cannot be a Heegaard surface.

Let $A_0$ be a horizontal annulus in $\Sigma'$.  The boundary of $A_0$ consists of two parallel simple closed curves in $F$.  Let $A_1$ be the component or the union of components of $\Sigma'$ adjacent to $A_0$.  The set $A_1$ is isotopic to the result of attaching a horizontal band to $\partial A_0 \times [0,1]$, or attaching a horizontal band from $\partial A_0 \times [0,1]$ to a third vertical annulus.  

If the horizontal band involves a third loop rather than a boundary component of $A_0$ then $A_0 \cup A_1$ is a three punctured sphere and (perhaps after an isotopy) is embedded in $F \times [0,1]$ such that the projection into $F$ is one-to-one.  In other words, we can think of $A_0 \cup A_1$ as the graph of a (Morse) function from a pair of pants to $\mathbf{R}$.  This function contains one index-two critical point (in $A_0$) and two index-one critical points (one in $A_0$ and one in $A_1$.)  There is a function on a pair of pants which agrees with this function on the boundary, but contains only one index-one critical point.  

Isotoping $A_0 \cup A_1$ to a graph of this new function, as in Figure~\ref{annremfig}, eliminates the horizontal annulus and reduces the number of components of $\Sigma \setminus F$.  Because we assumed $\Sigma'$ is minimal, we can rule out this case.  
\begin{figure}[htb]
  \begin{center}
  \includegraphics[width=3.5in]{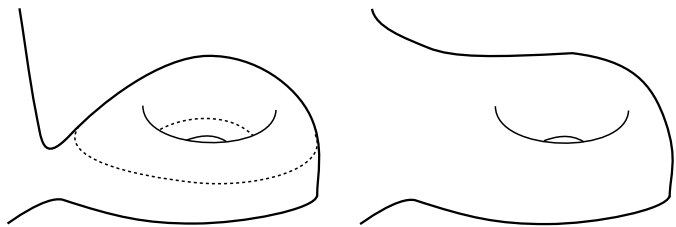}
  \caption{If a horizontal annulus is adjacent to an essential pair of pants which attaches the horizontal annulus to a separate loop, then the annulus and pair of pants can be isotoped into a single essential pair of pants.}
  \label{annremfig}
  \end{center}
\end{figure}

The horizontal band must go from a boundary loop in $A_0$ either back to itself or to the other loop. In the first case, the boundary of $A_1$ opposite from $\partial A_0$ consists of three (essential) simple closed curves which bound a three punctured sphere in $F$.  No two of these loops are parallel.  In the second case, the boundary of $A_1$ opposite $A_0$ is a single loop which bounds a onece punctured torus in $F$.

The surface $A_2$ adjacent to $A_1$, opposite $A_0$, consists of either vertical annuli or an essential pair of pants and possibly vertical annuli.  Because there is no pair of parallel loops in $\partial A_1$, the surface $A_2$ cannot contain a horizontal annulus.  If $A_2$ consists entirely of vertical annuli then we can push the saddle in $A_2$ through $F$ without changing the properties of $\Sigma'$ that we required at the beginning of the proof.  We can then push $A_0$ through $F$, reducing the number of components of $\Sigma'$.  Thus the minimality assumption implies that $A_2$ contains an essential pair of pants.

Let $A_3$ be the component or components of $\Sigma'$ adjacent to $A_2$ on the side opposite $A_1$.  Any horizontal annulus in $\Sigma'$ that is disjoint from $A_0 \cup A_1 \cup A_2$ must be a component of $A_3$ because we showed that such an annulus is adjacent to a pair of pants.  At least one component of $A_3$ must be a horizontal annulus in order for $\Sigma$ to be a closed surface.  Applying the same argument to $A_3$ and $A_2$ as we applied to $A_0$ and $A_1$ implies that $A_3$ is a single horizontal annulus.  The boundary opposite $A_3$ is either three loops that bound a pair of pants in $F$ or one loop that bounds a punctured torus.

In order for $\Sigma$ to close up, the boundary of $A_1$ opposite $A_0$ must have the same number of loops as the boundary of $A_2$ opposite $A_3$.  If each consists of one loop (which must be the same loop in $F$) or both consist of the same three loops in $F$ then $\Sigma$ is the union $A_0 \cup A_1 \cup A_2 \cup A_3$.  This surface is disjoint from a level surface of $\pi$, so there is an incompressible surface in the complement of $\Sigma$.  This contradicts the assumption that $\Sigma$ is a Heegaard surface.

Assume both boundaries consist of three loops each.  Because $A_1$ and $A_2$ are adjacent in $\Sigma$, at least one of these loops must be commong to both.  If one loop is common to both then $A_0 \cup A_1 \cup A_2 \cup A_3$ is a four punctured sphere.  If two loops are common to both then $A_0 \cup A_1 \cup A_2 \cup A_3$ is a twice punctured torus.  In either case, the surface is isotopic to a graph of a function on a four punctured sphere or a twice punctured torus.  On either surface, there is a function which has exactly two index-two critical points and the appropriate values on the boundary.  Isotoping $A_0 \cup A_1 \cup A_2 \cup A_3$ to the graph of such a function reduces the number of components of $\Sigma'$.  We assumed that the number of components is minimal, so we can rule out both cases.  We conclude that if $\Sigma$ is a Heegaard surface then $\Sigma$ is essentially embedded.
\end{proof}

\begin{Lem}
\label{essemblem}
If $\Sigma$ is an essentially embedded, separating, genus two surface in a surface bundle $M$ then the bundle structure on $M$ induces an almost bundle structure on each component of the complement of $\Sigma$.
\end{Lem}

\begin{proof}
Let $\Sigma$ be a genus two, essentially embedded surface in a surface bundle $M$.  Every regular level of the restriction map $\pi|_\Sigma$ is an essential loop in $\Sigma$.  There are no index 0 or 2 critical points because the level sets near a central singularity are all trivial loops.  The critical points in $\pi |_{\Sigma}$ must all be saddles.   The Euler characteristic of a genus two surface is $-2$ so there are exactly two saddle singularities.  

Let $x,y \in S^1$ and let $F = \pi^{-1}(\{x,y\})$.  By choosing $x$ and $y$ appropriately, we can ensure that each critical point of $\pi |_{\Sigma}$ is in a separate component of $M \setminus F$.  Then each component of $\Sigma \setminus F$ is either a vertical annulus or an essential pair of pants.

Let $H$ be the closure of a component of $M \setminus \Sigma$.  Then each component of $\Sigma \cap C$ is either an annulus in $\Sigma \setminus F$ or an essential pair of pants.  Exactly two components of $\partial H \setminus F$ contain one essential pair of pants each.  Let $C_0$ be the closure one of these component.  

There is a connected subsurface $F'$ of $F$ such that $C_0$ is homeomorphic to the complement in $F' \times [0,1]$ of a regular neighborhood $N$ of an essential arc in $F' \times\{0\}$ or $F' \times \{1\}$.    Without loss of generality, assume the arc is in $F' \times \{1\}$ and let $F''$ be the complement in $F' \times \{1\}$ of $N$.  (The surface may be connected or have two components.)  The complement $\partial C_0 \setminus \Sigma$ consists of  $F' \times \{0\}$ and $F''$.

Let $C'_0$ be the closure of a component of $H \setminus F$ such that $C_0 \cap C'_0$ is non-empty and assume $\partial C'_0 \cap \Sigma$ consists entirely of annuli.  The intersection $C'_0 \cap C_0$ in $C_0$ is either $F' \times \{0\}$ or a component of $F''$.  The bundle structure on $C'_0$ is that of a surface (with boundary) cross an interval, so $C_0 \cup C'_0$ is also homeomorphic to the complement in $F' \times [0,1]$ of $N$.  Define $C_1  = C_0 \cup C'_0$.

If there is a component $C'_1$ adjacent to $C_1$ such that $\partial C'_1 \cap \partial H$ consists of annuli and $C'_1 \cap C_1$ is connected then $C'_1 \cup C_1$ will be homeomorphic to $C_1$ and we define $C_2 = C_1 \cup C'_1$.  Continue in this fashion until we have constructed a subset $C_n$ such that for every component $C'_n$ adjacent to $C_n$, either the intersection $C'_n \cap C_n$ is not connected or  $\partial C'_n \cap \Sigma$ contains a pair of pants.  

If $\partial C'_n \cap \Sigma$ contains only annuli then $C'_n$ is homeomorphic to a surface cross an interval.  No component of $F''$ is homeomorphic to $F'$ so $C'_n \cap C_n$ cannot contain both $F' \times \{0\}$ and a component of $F''$.  In this case,  $C'_n \cap C_n$ must be connected. We conclude that $\partial C'_n \cap \partial H$ contains a pair of pants.  Because there is only one such component other than $C_0$ and $H$ is connected, we must have $C_n \cup C'_n = H$.  

By induction, $C_n$ is homeomorphic to $F' \times [0,1] \setminus N$.  Likewise, $C'_n$ is homeomorphic to a subsurface of $F$ cross an interval with a neighborhood of an arc removed.  Because $\partial C'_n \setminus \Sigma$ is homeomorphic to $\partial C_n \setminus \Sigma$, $C'_n$ must in fact be homeomorphic to $F' \times [0,1] \setminus N'$ where $N'$ is a regular neighborhood of an arc in $F' \times \{0\}$ and $(F' \times \{0\}) \setminus N'$ is homeomorphic to $F''$.

Let $X = C_n \cup_{F'} C'_n$ be the result of gluing the two components along $F' \times \{0\}$ in $C_n$ and $F' \times \{1\}$ in $C'_n$.  Then $X$ is homeomorphic to the complement in $F' \times [0,1]$ of a regular neighborhood of an arc $\alpha_0$ in $F' \times \{0\}$ and an arc $\alpha_1 \subset F' \times \{1\}$.  To complete the reconstruction of $H$, we simply glue $F' \times \{0\}$ to $F' \times \{1\}$ by some homeomorphism on the complements of $\alpha_0$ and $\alpha_1$, respectively.  Thus $H$ has the structure of an almost bundle.
\end{proof}

\begin{proof}[Proof of Lemma~\ref{bundlehomeolem}]
Let $(\Sigma, H_1, H_2)$ be a genus two Heegaard splitting for a surface bundle $M$ which is not a torus bundle or $S^1 \times S^2$.  This splitting is irreducible because $M$ is irreducible and does not allow a genus-one Heegaard splitting.  A weakly reducible, genus two Heegaard splitting is always reducible so $\Sigma$ is in fact strongly irreducible.   Thus Lemma~\ref{splitlem} implies that $\Sigma$ can be isotoped to an almost fair surface.   

By Lemma~\ref{almstfairlem}, $\Sigma$ must in fact be fair so Lemma~\ref{fairlem} implies that $\Sigma$ is isotopic to an essentially embedded surface.  Finally, Lemma~\ref{essemblem} implies that after this isotopy, the bundle structure of $M$ implies an almost bundle structure for either handlebody in the Heegaard splitting.  This completes the proof.
\end{proof}

\bibliographystyle{abbrv}
\bibliography{bundles}

\begin{thebibliography}{1}

\bibitem{bsc:bndls}
D.~Bachman and S.~Schleimer.
\newblock {Surface bundles versus Heegaard splittings}.
\newblock {\em Communications in Analysis and Geometry 13}, 5:1--26, 2005.

\bibitem{sch:solv}
D.~Cooper and M.~Scharlemann.
\newblock {The structure of a solvmanifold's Heegaard splittings}.
\newblock {\em Turkish J. Math}, 23(1):1--18, 1999.

\bibitem{Kob:disc}
T.~Kobayashi and O.~Saeki.
\newblock {The Rubinstein-Scharlemann graphic of a 3-manifold as the
  discriminant set of a stable map}.
\newblock {\em Pacific Journal of Mathematics}, 195(1):101--156, 2000.

\bibitem{lum:niels}
M.~Lustig and Y.~Moriah.
\newblock {Nielsen equivalence in Fuchsian groups and Seifert fibered spaces}.
\newblock {\em Topology}, 30(2):191--204, 1991.

\bibitem{rub:compar}
H.~Rubinstein and M.~Scharlemann.
\newblock {Comparing Heegaard Splittings of Non-Haken 3-Manifolds}.
\newblock {\em Topology}, 35(4):1005--1026, 1996.

\bibitem{sch:detect}
M.~Scharlemann.
\newblock {Local detection of strongly irreducible Heegaard splittings}.
\newblock {\em Topology Appl.}, 90(1-3):135--147, 1998.

\bibitem{thr:norm}
W.~Thurston.
\newblock A norm for the homology of 3-manifolds.
\newblock {\em Mem. Amer. Math. Soc.}, 59(339):99--130, 1986.

\end{thebibliography}

\end{document}